\documentclass{elsart}
\usepackage{amssymb}
\usepackage{amsmath}
\usepackage{graphicx}
\usepackage{epsfig}
\usepackage[francais]{babel}
\usepackage[latin1]{inputenc}
\begin{document}

\begin{frontmatter}
\title{Un Mod\`{e}le simple d'injection diphasique avec phase condensable}


\author[lab1]{G\'erard Debicki}
\author[lab2]{J\'er\^ome Pousin}
\author[lab2]{Eric Zeltz}
\corauth[cor1]{Corresponding Author Fax: 00 33 4 72 43 85 29}

\address[lab1]{ Universit\'e de Lyon CNRS  \\INSA-Lyon URGC
}
\address[lab2]{Universit\'e de Lyon CNRS  \\INSA-Lyon ICJ UMR 5208, bat. L. de
Vinci,
\\20 Av. A. Einstein, F-69100 Villeurbanne Cedex France }

\begin{abstract}
L'objectif de cette note est de proposer un modèle mathématique simple permettant de comprendre l'arrêt de la pénétration d'un flux de vapeur d'eau condensable sur un mur de béton, observé expériementalement (voir par exemple  les diff\'{e}rentes exp\'{e}riences d'injection de vapeur dans du b\'{e}ton présentées dans  \cite{sh}, et \cite{la}). Un modèle homogénéisé simple d'injection dans un milieu poreux est proposé, donnant une borne pour la position asymptotique en temps du front de pénétration.

\noindent The aim of this paper is to propose a simple mathematical model to understand the decision of the penetration of a stream of water vapor condensing on a concrete wall, observed experimentally (see for example the situations described in \cite{sh}, et \cite{la}). A  simple homogenized model for the injection in a porous medium is proposed, giving a bound for the asymptotic-time position at the front of penetration.
\end{abstract}
\end{frontmatter}
\section{Introduction}

 Un mod\`{e}le simple d'injection
d'un fluide condensable incompressible dans un milieu poreux empli d'un
autre fluide incompressible est présenté dans \cite{pous}. Pour cela les auteurs ont adapt\'{e} \`{a}
cette situation un autre mod\`{e}le (cf. \cite{maa}) qu'ils avaient obtenu
par homog\'{e}n\'{e}isation \`{a} partir du syst\`{e}me de Navier-Stokes et
qui d\'{e}crivait l'injection isotherme d'un fluide incompressible non
condensable dans un milieu poreux empli d'un autre fluide incompressible.

L'adaptation consistait \`{a} d\'{e}finir la vitesse du fluide condensable
entrant comme somme de la vitesse qu'il aurait s'il n'\'{e}tait pas
condensable ( donc v\'{e}rifiant la mod\'{e}lisation d\'{e}crite dans \cite%
{maa}) et d'un terme de "recul" traduisant la diminution de volume provoqu%
\'{e}e par la condensation.

Cela leur a permis notamment de montrer que la vitesse du d\'{e}placement de
l'interface entre les deux fluides \'{e}tait d\'{e}croissante et
asymptotiquement nulle pour les temps infinis.

Ainsi ils purent avancer une explication pour expliquer l'arr\^{e}t observ%
\'{e} exp\'{e}rimentalement de l'avanc\'{e}e d'une vapeur condensable inject%
\'{e}e dans du b\'{e}ton. Explication qui compl\`{e}te celle qui \'{e}tait
jusqu'alors donn\'{e}e dans la litt\'{e}rature pour expliquer ce ph\'{e}nom%
\`{e}ne: \`{a} savoir, la formation d'un "bouchon" de condensation (voir par
exemple \cite{con},\cite{sh} ou \cite{la}).

Tout cela est par ailleurs d\'{e}velopp\'{e} et d\'{e}taill\'{e} dans la th%
\`{e}se de Zeltz (cf. \cite{zel}) .

La d\'{e}marche adopt\'{e}e par ces auteurs, bri\`{e}vement r\'{e}sum\'{e}e,
est la suivante:

Ils consid\'{e}rent que la condensation provient essentiellement de la diff%
\'{e}rence entre la pression r\'{e}elle et la pression de vapeur saturante
du fluide.

Et que donc cette condensation ne se fait pas uniquement au niveau de
l'interface entre le fluide inject\'{e} condensable et le fluide r\'{e}%
siduel non condensable. Mais qu'elle se fait de fa\c{c}on homog\`{e}ne dans
toute la portion poreuse mouill\'{e}e par le fluide inject\'{e}.

Plus pr\'{e}cis\'{e}ment, le terme utilis\'{e} par ces auteurs pour mod\'{e}%
liser ce recul est de la forme lin\'{e}aire $\delta V(t)$, $V(t)$ \'{e}tant
le volume du fluide entr\'{e} \`{a} l'instant $t$ et $\delta $ un
coefficient li\'{e} \`{a} la capacit\'{e} de condensation du fluide consid%
\'{e}r\'{e}, \`{a} la temp\'{e}rature donn\'{e}e.

L'id\'{e}e pr\'{e}cise sur laquelle repose ce choix peut s'exprimer ainsi:

Si $p_{i}$ repr\'{e}sente la pression appliqu\'{e}e \`{a} la vapeur inject%
\'{e}e et $p_{vs}$ la pression de vapeur saturante, la diff\'{e}rence $%
p_{i}- $ $p_{vs}$ est constante puisqu'en condition isotherme $p_{vs}$ est
constante. Or la condensation est proportionnelle \`{a}
\begin{tabular}{l}
$p_{i}-$ $p_{vs}$%
\end{tabular}
(cf. par exemple \cite{hu}). Et donc dans ces conditions, la quantit\'{e} de
fluide qui se condense dans le milieu poreux est effectivement et uniquement
proportionnelle au volume du fluide d\'{e}j\`{a} entr\'{e}.

La question \'{e}tudi\'{e}e dans les paragraphes qui suivent est celle-ci:

En gardant le point de vue qui vient d'\^{e}tre d\'{e}crit, le r\'{e}sultat
obtenu dans \cite{pous} (d\'{e}c\'{e}l\'{e}ration progressive de l'avanc\'{e}%
e du fluide condensable inject\'{e} jusqu'\`{a} une vitesse nulle ) est-il g%
\'{e}n\'{e}ralisable dans le cas d'une injection non isotherme?

\subsection{Description du milieu poreux et des autres hypoth\`{e}ses}

\bigskip

Le but principal de cette \'{e}tude est d'obtenir une mod\'{e}lisation \`{a}
la fois simple et fiable d'une injection de vapeur dans une enceinte en b%
\'{e}ton empli d'air humide.

La structure poreuse envisag\'{e}e et les hypoth\`{e}ses physiques ont donc
\'{e}t\'{e} choisies en fonction.

\begin{figure}[h]
\begin{center}
\begin{minipage}{0.48\linewidth}
\epsfig{file=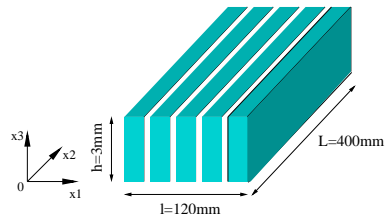, width=5cm}
\end{minipage}
\caption{Sch\'{e}matisation du milieu poreux }\label{pores}
\end{center}
\end{figure}

La mod\'{e}lisation de la structure poreuse est celle qui a \'{e}t\'{e}
utilis\'{e}e dans \cite{pous} et dans \cite{maa}:

Les pores ont la forme de "slits", c'est-\`{a}-dire de fentes tr\`{e}s
fines. Ils sont suppos\'{e}s parall\`{e}les et non interconnect\'{e}s. La
figure \ref{pores} sch\'{e}matise un tel milieu.

\bigskip

Chaque slit est donc un parall\'{e}lipip\`{e}de de dimensions $e\times
h\times L$ de l'espace $(O,\overrightarrow{\mathbf{x}_{1}},\overrightarrow{%
\mathbf{x}_{2}},\overrightarrow{\mathbf{x}_{3}})$

avec $e<<h<<L.$

L'\'{e}coulement dans ces slits peut alors \^{e}tre assimil\'{e} \`{a} un
\'{e}coulement de Hele-Shaw, c'est-\`{a}-dire un \'{e}coulement visqueux
entre deux plans parall\`{e}les voisins.

Ce qui permet, en proc\'{e}dant de la m\^{e}me mani\`{e}re que dans \cite%
{ric}, de se ramener \`{a} des pores 2D $(x_{2},x_{3})\in \Omega
=]0,L[\times ]0,h[$ du plan $(O,\overrightarrow{\mathbf{x}_{2}},%
\overrightarrow{\mathbf{x}_{3}}).$

\begin{figure}[h]
\begin{center}
\epsfig{file=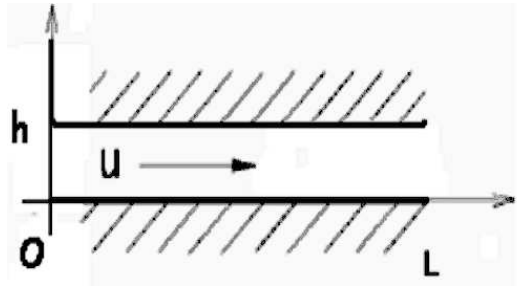, width=10cm}
\caption{Sch\'{e}ma d'un pore.}\label{unpore}
\end{center}
\end{figure}

Un tel pore est sch\'{e}matis\'{e} dans la figure \ref{unpore}.

Le milieu poreux, et donc aussi chacun des pores, est initialement rempli
d'un premier fluide r\'{e}siduel $f_{r\text{ }}$\`{a} pressions et temp\'{e}%
ratures ambiantes. Ce fluide est suppos\'{e} newtonien, visqueux et
incompressible. Cela peut \^{e}tre consid\'{e}r\'{e} comme \'{e}tant le cas
pour de l'air en d\'{e}placement lent, car alors sa compressibilit\'{e} est n%
\'{e}gligeable (pour la validit\'{e} de cette approximation, on pourra voir
par exemple \cite{boi}).

Le fluide inject\'{e} $f_{i\text{ }}$est aussi newtonien, visqueux et
incompressible. Il a la propri\'{e}t\'{e} suppl\'{e}mentaire qu'il se
condense lorsque sa pression d\'{e}passe la pression de vapeur saturante.
C'est le cas de la vapeur d'eau. Ses pressions et temp\'{e}ratures \`{a}
l'injection sont constantes et suppos\'{e}es tr\`{e}s sup\'{e}rieures aux
pressions et temp\'{e}ratures ambiantes.

Dans la suite, la notation suivante sera syst\'{e}matiquement utilis\'{e}e:
ce qui se rapporte au fluide r\'{e}siduel $f_{r\text{ }}$ sera indic\'{e}
par $_{r\text{ }}$ et ce qui se rapporte au fluide inject\'{e} $f_{i}$ par $%
_{i}$.

Pour les viscosit\'{e}s $\eta _{i}$ et $\eta _{r}$ des deux fluides en jeu,
elles sont assez voisines pour pouvoir les confondre, ce qui est le cas en
particulier pour l'air humide et la vapeur d'eau (cf. \cite{in}). Donc ces
viscosit\'{e}s sont suppos\'{e}es v\'{e}rifier:

\begin{tabular}{l}
$\eta _{i}(x_{2},x_{3},t)=\eta _{r}(x_{2},x_{3},t)=:\eta (x_{2},x_{3},t)%
\mbox{ dans }\;\Omega \times ]0,T[$%
\end{tabular}%
.

Puisque ceux-ci sont essentiellement sous phase gazeuse, les tensions
superficielles sont n\'{e}gligeables et il n'en est donc pas tenu compte.

De m\^{e}me, les forces de gravitation peuvent \^{e}tre consid\'{e}r\'{e}es
comme minimes pour des gaz \'{e}voluant dans un tel milieu poreux, et sont
donc n\'{e}glig\'{e}es.

Les conditions d'immiscibilit\'{e} des deux fluides en jeu, de transmission
\`{a} l'interface (continuit\'{e}s de la vitesse et de la composante normale
du tenseur de contrainte) et d'adh\'{e}rence aux parois (vitesse nulle au
bord) sont celles qui sont habituellement retenues dans ce genre d'\'{e}%
coulement de fluides newtoniens visqueux incompressibles.

Elles sont identiques \`{a} celles d\'{e}j\`{a} utilis\'{e}es dans \cite%
{pous}.

Par contre, et c'est la seule diff\'{e}rence dans les hypoth\`{e}ses avec
celles adopt\'{e}es dans \cite{pous}, les conditions ne sont pas isothermes
et l'\'{e}volution de la temp\'{e}rature provoqu\'{e}e par la conduction de
la chaleur \`{a} travers le milieu poreux est prise en compte. Cela sera pr%
\'{e}cis\'{e} par la suite.

\subsection{\protect\bigskip D\'{e}veloppement asymptotique et caract\'{e}%
risation de la vitesse}

A tout instant $t\geq 0$, le domaine $\;\Omega $ est r\'{e}union disjointe
de la partie $\;\Omega _{r}(t)$ encore emplie du fluide non condensable r%
\'{e}siduel et de la partie $\Omega _{i}(t)$ d\'{e}j\`{a} emplie du fluide
entrant condensable:
\begin{tabular}{l}
$\Omega =\Omega _{r}(t)\cup \Omega _{i}(t)$%
\end{tabular}
avec
\begin{tabular}{l}
$\Omega _{r}(t)\cap \Omega _{i}(t)=\varnothing $%
\end{tabular}%
$.$

L'interface $\Gamma (t)$ \`{a} l'instant $t$ correspond donc \`{a}
\begin{tabular}{l}
$\overline{\Omega _{r}(t)}\cap \overline{\Omega _{i}(t)}$%
\end{tabular}%
$.$

\begin{figure}[h]
\begin{center}
\epsfig{file=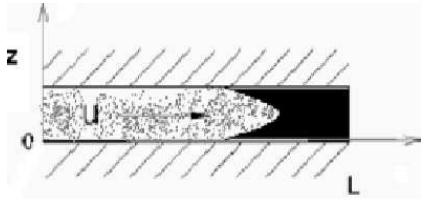, width=10cm}
\caption{Injection de la vapeur dans le pore. }\label{condensation}
\end{center}
\end{figure}

Cela est sch\'{e}matis\'{e} dans la figure \ref{condensation}.

La vitesse ${\vec{v}}$ dans $\Omega $ peut \^{e}tre d\'{e}compos\'{e}e par:%
\newline
\qquad \qquad\ ${\vec{v}(x}_{2},z,t)$ ${=}$ $\overrightarrow{{w}}{(x}%
_{2},z,t)-k(x_{2},t)$ $\overrightarrow{\mathbf{x}_{2}}$

avec $\overrightarrow{{w}}$ et $k$ ainsi d\'{e}finis:

\begin{itemize}
\item $\overrightarrow{{w}}$ repr\'{e}sente la vitesse qu'on aurait pour le m%
\^{e}me probl\`{e}me mais sans condensation.

\item $k(x_{2},t)$ est un r\'{e}el positif ou nul qui correspond \`{a} la
perte de vitesse provoqu\'{e}e par la diminution de volume qu'entraine la
condensation.
\end{itemize}

Les d\'{e}veloppements asymptotiques sont \'{e}crits en fonction de $%
\varepsilon =h<<L$\newline

\`{a} partir du domaine
\begin{tabular}{l}
$\Omega ^{\varepsilon }=]0,L[\times ]0,\varepsilon \lbrack $%
\end{tabular}
.

Les inconnues que constituent la vitesse ${\vec{v}}^{\varepsilon }$ du
fluide, sa pression $p^{\varepsilon }$ et sa viscosit\'{e} $\eta
^{\varepsilon }$ sont red\'{e}finies sur $\Omega =\Omega ^{1}$ et d\'{e}%
pendent des variables $x_{2},$ $z=\frac{x_{3}}{\varepsilon }$ et $t$, comme
cela est fait classiquement dans la m\'{e}thode de r\'{e}duction introduite
par Ciarlet et \textit{al.} (cf. \cite{cia}).

Les d\'{e}veloppements utilis\'{e}s sont pr\'{e}cis\'{e}ment ceux-ci:\medskip

$\left\{
\begin{tabular}{lll}
$\eta ^{\varepsilon }(x_{2},z,t)$ & $=$ & $\underset{k\geq 0}{\sum }%
\varepsilon ^{k}\eta ^{k}(x_{2},z,t)$ \\
&  &  \\
$p^{\varepsilon }(x_{2},z,t)$ & $=$ & $\varepsilon ^{-2}p^{0}(x_{2},t)+%
\underset{k\geq 0}{\sum }\varepsilon ^{k-2}p^{k}(x_{2},z,t)$ \\
&  &  \\
${\vec{v}}^{\varepsilon }(x_{2},z,t)$ & $=$ & $\underset{k\geq 0}{\sum }%
\varepsilon ^{k}{\vec{v}}^{k}(x_{2},z,t)$ $\mathit{avec}$ ${\vec{v}}^{k}({v}%
_{2}^{k},v_{z}^{k})$%
\end{tabular}%
\ \right. \medskip \medskip $

\noindent

A partir du syst\`{e}me de Stokes et comme cela est montr\'{e} dans le
chapitre 4 de \cite{zel}, le syst\`{e}me homog\'{e}n\'{e}is\'{e} obtenu en
consid\'{e}rant uniquement les termes d'ordre z\'{e}ro est alors le suivant:

\begin{equation}
\left\{
\begin{array}{lll}
-\dfrac{\partial }{\partial z}\{\eta ^{0}\dfrac{\partial v_{2}^{0}}{\partial
z}\}+\dfrac{\partial p^{0}}{\partial x_{2}} & =0 & \mbox{ dans }\;\Omega
\times ]0,T[ \\
\left( \dfrac{\partial v_{z}^{1}}{\partial z}+\dfrac{\partial v_{2}^{0}}{%
\partial x_{2}}\right) & =-\dfrac{\partial }{\partial x_{2}}\left\{
k^{0}(x_{2},t)\right\} & \mbox{ dans }\;\Omega \times ]0,T[ \\
\dfrac{\partial \eta ^{0}}{\partial t}+\dfrac{\partial }{\partial z}\{\eta
^{0}v_{z}^{1}\}+\dfrac{\partial }{\partial x_{2}}\{\eta ^{0}v_{2}^{0}\} & =-%
\dfrac{\partial }{\partial x_{2}}\left\{ \eta ^{0}k^{0}(x_{2},t)\right\} & %
\mbox{ dans }\;\Omega \times ]0,T[%
\end{array}%
\right.  \label{hom}
\end{equation}%
\medskip avec pour conditions initiales et au bord:\medskip\
\begin{equation}
\left\{
\begin{array}{l}
\vec{v}^{0}(0,z,t)=(u_{00},0)\;\mbox{ dans }]0,1[\times ]0,T[ \\
\vec{v}^{0}(L,z,t)=(u_{0L},0)\;\mbox{ dans }]0,1[\times ]0,T[ \\
\int_{0}^{1}u_{00}\ dz=q\;\mbox{ sur }\;]0,T[ \\
\vec{v}^{0}(x_{2},z,t)=-k^{0}(x_{2},t)\overrightarrow{\mathbf{x}_{2}}\;%
\mbox{ dans }]0,L[\times \left\{ 0,1\right\} \times ]0,T[ \\
\eta ^{0}(x_{2},z,0)=\eta _{r}\mbox{ sur  }]0,L[\times ]0,1[ \\
\eta ^{0}(0,z,t)=\eta _{i}\;\mbox{ sur }]0,1[\times ]0,T[%
\end{array}%
\right.  \label{bord}
\end{equation}

\subsection{Comportement asymptotique de l'interface}

L'obtention du comportement asymptotique de l'interface se fait en supposant que la température est donnée par la solution de l'équation de la chaleur modélisant l'\'{e}volution de la temp\'{e}rature provoqu%
\'{e}e par conduction en milieu unidimensionnel semi infini.
Nous préciserons d'abord  les conditions initiales et au bord
de la temp\'{e}rature et de la pression.
L'\'{e}volution de la pression de vapeur saturante en fonction de la
temp\'{e}rature est alors d\'{e}duite.
Ensuite est d\'{e}termin\'{e}e l'\'{e}volution de pression soumise \`{a}
la vapeur \`{a} l'int\'{e}rieur du milieu poreux.
Nous définissons alors la  fonction de condensation $\delta
:(x_{2};t)\longmapsto $ $\delta (x_{2};t).$ Cela permet de
d\'{e}terminer la fonction $k(x_{2},t)$ intervenant dans la définition de la vitesse.

Finalement nous  pr\'{e}cisons la position $x_{2}(z,t)$ de
l'interface et nous \'{e}tudions son comportement asymptotique quand $t$ tend
vers $+\infty .$

\subsubsection{Conditions initiales et au bord}

A l'entr\'{e}e $E=\{x_{2}=0\}$, les temp\'{e}ratures $\Theta $ et pressions $%
p$ sont maintenues \`{a} partir de l'instant initial $t=0$ aux valeurs
constantes $\Theta (0,t):=\Theta _{E}$ et $p(0,t):=p_{E}.$ De m\^{e}me, apr%
\`{e}s la sortie $S=\{x_{2}=L\}$, la pression est maintenue \`{a} la valeur
constante $p(L,t):=p_{S}.$

Enfin, \`{a} l'int\'{e}rieur du milieu poreux, les conditions initiales de
temp\'{e}rature et de pression sont les suivantes: Pour tout $x_{2}\in ]0;L]:%
\begin{tabular}{l}
$\Theta (x_{2},0)=\Theta _{S}$ et $p(x_{2},0)=p_{S}$%
\end{tabular}%
.$

En raison des probl\`{e}mes concrets qui sont mod\'{e}lis\'{e}s ici, ces
valeurs v\'{e}rifient: $\Theta _{E}>\Theta _{S}$ et $p_{E}>p_{S}.$

\subsubsection{D\'{e}termination de la temp\'{e}rature et de sa limite
asymptotique}

Les transferts de chaleur dans le milieu poreux consid\'{e}r\'{e} sont suppos%
\'{e}s se faire uniquement par conduction thermique et non par convection.
Cette approximation est valide si la p\'{e}n\'{e}tration du fluide dans le
milieu poreux est tr\`{e}s lente et si les temp\'{e}ratures du fluide inject%
\'{e} d\'{e}passent 100$%
{{}^\circ}%
C$ comme cela est clairement \'{e}tabli dans l'\'{e}tude param\'{e}trique r%
\'{e}sum\'{e}e pages 195 et suivantes dans la th\`{e}se de Laghcha (cf. \cite%
{la}). D'autres justifications de la validit\'{e} de cette approximation se
trouvent par exemple dans \cite{ric} ou encore dans \cite{bi}.

La temp\'{e}rature $\Theta $ est donc suppos\'{e}e ob\'{e}ir uniquement \`{a}
une \'{e}quation de chaleur issue de la loi de Fourier.

Cette loi est suppos\'{e}e de plus unidirectionnelle, puisque les flux de
fluides dans les exp\'{e}riences consid\'{e}r\'{e}es sont orient\'{e}s des
masses chaudes de l'entr\'{e}e vers les masses froides de la sortie.

Elle est donc du type $\frac{\partial ^{2}\Theta }{\partial x_{2}^{2}}=\frac{%
1}{K}\frac{\partial \Theta }{\partial t}$ o\`{u} $K$ est une constante
strictement positive d\'{e}pendant uniquement des propri\'{e}t\'{e}s
physiques du milieu poreux consid\'{e}r\'{e}, en l'occurence la densit\'{e},
la conductivit\'{e} et la chaleur sp\'{e}cifique du mat\'{e}riau solide qui
le constitue (cf. par exemple \cite{in}).

Compte-tenu de ce qui pr\'{e}c\`{e}de, la temp\'{e}rature $\Theta $ est
solution du probl\`{e}me suivant:

\begin{equation}
\left\{
\begin{tabular}{lll}
$\frac{\partial ^{2}\Theta }{\partial x_{2}^{2}}$ & $=\frac{1}{K}\frac{%
\partial \Theta }{\partial t}$ & $\mbox{ dans }\;0< x_2; \,   0<t<T$; \\
$\Theta (0,t)$ & $=\Theta _{E}>\Theta _{S}$ & pour $0<t<T$; \\
$\Theta (x_{2},0)$ & $=\Theta _{S}$ & pour $0<x_{2}$
\end{tabular}%
\begin{tabular}{lll}
&  &  \\
&  &  \\
&  &
\end{tabular}%
\right.   \label{temper}
\end{equation}%
\medskip

Ce probl\`{e}me est classique et admet pour solution fondamentale (cf. \cite%
{boi} page 106):

$\
\begin{tabular}{l}
$\Theta (x_{2},t)=\Theta _{S}+(\Theta _{E}-\Theta _{S})\left( 1-\frac{2}{%
\sqrt{\pi }}\int_{0}^{\frac{x_{2}}{2\sqrt{Kt}}}e^(-x^{2})dx\right) $%
\end{tabular}%
\medskip $

et impose donc la valeur suivante \`{a} $\Theta (L,t):$

$\Theta (L,t)=\Theta _{S}+(\Theta _{E}-\Theta _{S})\left( 1-\frac{2}{\sqrt{%
\pi }}\int_{0}^{\frac{L}{2\sqrt{Kt}}}e^(-x^{2})dx\right) $

Remarque 1: \textit{Pour tout} $x_{2}\in \lbrack 0;L]:\underset{t\rightarrow
+\infty }{\lim }\Theta (x_{2},t)=\Theta _{E}$ \textit{et pour tout }$t$%
\textit{\ fix\'{e} dans} $%
\mathbb{R}
^{+}:\underset{L\rightarrow +\infty }{\lim }\Theta (L,t)=\Theta _{S}.$

Remarque 2: \textit{La viscosit\'{e} pour un fluide donn\'{e} n'\'{e}tant
fonction que de la temp\'{e}rature, et celle-ci \'{e}tant ici fonction
uniquement de }$x_{2}$\textit{\ et de }$t$\textit{, la viscosit\'{e} commune
des deux fluides v\'{e}rifie donc: }

\begin{tabular}{l}
$\eta (x_{2},z,t)=\eta (x_{2},t)\mbox{ dans }\;\Omega \times ]0,T[$%
\end{tabular}%
\medskip

\subsubsection{Evolution de la pression de vapeur saturante et de sa limite
asymptotique}

En supposant que la vapeur en jeu se comporte comme un gaz parfait, la
formule de Clapeyron (cf. encore \cite{in}) permet d'\'{e}valuer la pression
de vapeur saturante par:
\begin{tabular}{l}
$p_{vs}=p_{0}e^{\lambda (\frac{1}{\Theta _{0}}-\frac{1}{\Theta })}$%
\end{tabular}

avec :

$\Theta _{0}$ : temp\'{e}rature d'\'{e}bullition de la substance \`{a} une
pression $p_{0}$ donn\'{e}e, en degr\'{e}s $Kelvin.$

$p_{vs}$ : pression de vapeur saturante, dans la m\^{e}me unit\'{e} que $%
p_{0}$

$\lambda :$ une constante strictement positive d\'{e}pendant de la masse
molaire et de la chaleur latente de la substance.

\bigskip En prenant pour temp\'{e}rature de r\'{e}f\'{e}rence $\Theta _{0}$
la temp\'{e}rature $\Theta _{S}$ et en notant $\pi _{S}$ la pression pour
laquelle le fluide $f_{i}$ est \`{a} \'{e}bullition \`{a} la temp\'{e}rature
$\Theta _{S}$, la pression de vapeur saturante du fluide $f_{i}$ v\'{e}rifie
donc
\begin{tabular}{l}
$p_{vs}=\pi _{S}e^{\lambda _{i}(\frac{1}{\Theta _{S}}-\frac{1}{\Theta })}$%
\end{tabular}
o\`{u} $\lambda _{i}$ est la constante $\lambda $ associ\'{e}e \`{a} ce
fluide $f_{i}.$

Compte-tenu de la Remarque 1, il vient:

Pour tout $x_{2}\in \lbrack 0;L]:%
\begin{tabular}{l}
$\underset{t\rightarrow +\infty }{\lim }p_{vs}(x_{2},t)=\pi _{S}e^{\lambda
_{i}(\frac{1}{\Theta _{S}}-\frac{1}{\Theta _{E}})}$%
\end{tabular}%
.$

\subsubsection{D\'{e}termination de la pression appliqu\'{e}e et de sa
limite asymptotique}

La  relation $%
\begin{array}{lll}
-\dfrac{\partial }{\partial z}\{\eta ^{0}\dfrac{\partial v_{2}^{0}}{\partial
z}\}+\dfrac{\partial p^{0}}{\partial x_{2}} & =0 & \mbox{ dans }\;\Omega
\times ]0,T[%
\end{array}%
$ dans le syst\`{e}me \ref{hom}, obtenue par homog\'{e}n\'{e}isation \`{a} des équations de Stokes (cf. \cite{pous} et \cite{zel} par exemple) permet de calculer la
pression $p_{i}$ impos\'{e}e \`{a} la vapeur (qui est dans la suite
confondue avec sa pression essentielle $p_{0},$ le terme d'ordre $0$ de son d%
\'{e}veloppement asymptotique).

$v_{2}^{0}$ vérifiant  le syst\`{e}me \ref{hom} est donnée par (voir proposition 15 p. 52 \cite{zel})
:\medskip

\begin{equation}
v_{2}^{0}(x_{2},z,t)=q\frac{P(x_{2},z)}{R(x_{2})}-k^{0}(x_{2},t)
\label{schnell}
\end{equation}%
\medskip

o\`{u} les fonctions $P$ et $R$ sont ainsi d\'{e}finies:

$\left\{
\begin{array}{lll}
P(x_{2},z,t) & = & (\int_{0}^{z}\frac{d\xi }{\eta ^{0}(x_{2},\xi ,t)}%
)(\int_{0}^{1}\frac{\xi d\xi }{\eta ^{0}(x_{2},\xi ,t)})-(\int_{0}^{z}\frac{%
\xi d\xi }{\eta ^{0}(x_{2},\xi ,t)})(\int_{0}^{1}\frac{d\xi }{\eta
^{0}(x_{2},\xi ,t)}) \\
&  & \; \\
R(x_{2},t) & = & (\int_{0}^{1}\frac{\xi ^{2}d\xi }{\eta ^{0}(x_{2},\xi ,t)}%
)(\int_{0}^{1}\frac{d\xi }{\eta ^{0}(x_{2},\xi ,t)})-(\int_{0}^{1}\frac{\xi
d\xi }{\eta ^{0}(x_{2},\xi ,t)})^{2}%
\end{array}%
\right. .\medskip $

En utilisant la Remarque 2, un calcul simple montre que:\newline
\begin{tabular}{lll}
$P(x_{2},z,t)=\frac{\eta ^{0}}{2}z-\frac{\eta ^{0}}{2}z%
{{}^2}%
$ & et & $R(x_{2},t)=\frac{\eta ^{0}}{2}$%
\end{tabular}
$.\medskip $

Ceci donne:%
\begin{equation}\label{vitesse}
v_{2}^{0}(x_{2},z,t)=q(z-z^{2})-k^{0}(x_{2},t)
\end{equation}

Ainsi nous évaluons la dérivée de la pression avec:%
\begin{tabular}{l}
$\dfrac{\partial p^{0}}{\partial x_{2}}=\dfrac{\partial }{\partial z}\{\eta
^{0}\dfrac{\partial v_{2}^{0}}{\partial z}\}=-2q\eta ^{0}(x_{2},t)$%
\end{tabular}%
$.\medskip $

Donc la pression de la vapeur inject\'{e}e v\'{e}rifie:

\begin{equation}
p_{i}(x_{2},t)=-2q\int_{0}^{x_{2}}\eta ^{0}(\zeta ,t)d\zeta +p_{E}
\label{pr}
\end{equation}

Il reste à préciser  la viscosit\'{e} $\eta (x_{2},t)= \eta
^{0}(x_{2},t)$ pour pr\'{e}ciser cette pression.

Or elle-m\^{e}me d\'{e}pend de la temp\'{e}rature suivant la loi de
Sutherland (cf. \cite{in}):\medskip

\begin{tabular}{l}
$\eta ^{0}(x_{2},t)=\eta _{\Theta _{E}}\left( \frac{\Theta (x_{2},t)}{\Theta
_{E}}\right) ^{1.5}\frac{\Theta _{E}+\Psi }{\Theta (x_{2},t)+\Psi }$%
\end{tabular}
\medskip \newline
o\`{u} l'on note $\eta _{\Theta _{E}}$ la viscosit\'{e} du fluide \`{a} la
temp\'{e}rature d'entr\'{e}e $\Theta _{E}$ et o\`{u} $\Psi $ est une
constante d\'{e}pendant du fluide visqueux consid\'{e}r\'{e}.

Il en r\'{e}sulte que pour tout $x_{2}\in \lbrack 0;L]:%
\begin{tabular}{l}
$\underset{t\rightarrow +\infty }{\lim }\eta ^{0}(x_{2},t)=\eta
^{0}(x_{2},\infty )=\eta _{\Theta _{E}}$%
\end{tabular}%
$

Mais par ailleurs, la pression $p_{i}$ de la vapeur inject\'{e}e d\'{e}croit
jusqu'\`{a} la valeur $p_{S}$ qui est impos\'{e}e par la sortie, et reste
donc ensuite constante \`{a} cette valeur.

Donc, en tenant compte de cette remarque et en passant \`{a} la limite dans %
\ref{pr}, la pression $p_{i}$ tend vers une limite asymptotique $%
p_{i}(x_{2},\infty ):=\underset{t\rightarrow +\infty }{\lim }p_{i}(x_{2},t)$
ainsi d\'{e}finie d'après \ref{pr}:

Pour tout $x_{2}\in \lbrack 0;L],$\newline

\begin{equation}
p_{i}(x_{2},\infty )=\left\{
\begin{tabular}{l}
$-2q\eta _{\Theta _{E}}x_{2}+p_{E}$ si $x_{2}\leq \frac{p_{E}-p_{S}}{2q\eta
_{\Theta _{E}}}$ \\
$p_{S}$ sinon.%
\end{tabular}%
\right.  \label{PRL}
\end{equation}

Remarque 3: \textit{La pression asymptotique }$p_{i}(x_{2},\infty )$ \textit{%
obtenue par cette mod\'{e}lisation est donc essentiellement une fonction
affine d\'{e}croissante de la variable d'espace} $x_{2}$. \textit{Cela est
en assez bonne correspondance avec le comportement de la pression observ\'{e}%
e dans les diff\'{e}rentes exp\'{e}riences d'injection de vapeur dans du b%
\'{e}ton faites par Shekarchi et par Laghcha dans leurs th\`{e}ses
respectives \cite{sh}, et \cite{la}. Voir par exemple les figures V.29 et
V.30 page 202 de \cite{sh} et les figures IV. 6 et IV. 7 pages 186-187 de
\cite{la}. }

\subsubsection{D\'{e}termination de la fonction de condensation $\protect%
\delta $ et de sa limite asymptotique}

\bigskip

Nous avons supposé que la condensation se faisait continuement jusqu'\`{a} l'avanc%
\'{e}e du front et non seulement au niveau de celui-ci.

En effet, le ph\'{e}nom\`{e}ne de condensation de la vapeur provient
essentiellement de l'existence d'une diff\'{e}rence positive entre la
pression impos\'{e}e au fluide et sa pression de vapeur saturante (cf. par
exemple \cite{hu}).

Et cet \'{e}cart est pr\'{e}sent sur toute la partie mouill\'{e}e par le
fluide condensable inject\'{e}, et pas seulement au niveau de l'interface
avec le fluide r\'{e}siduel.

Aussi le calcul de la fonction de condensation $\delta $ est le suivant:

Soit $p_{i}$ la pression initiale de la vapeur inject\'{e}e.

La vapeur qui se condense est celle qui est n\'{e}cessaire pour que $p_{i}$
prenne la valeur de la vapeur
saturante $p_{vs}(\Theta )$.

Plus pr\'{e}cis\'{e}ment:
\begin{tabular}{l}
$\delta (x_{2},t)=\max \left( 0;\frac{p_{i}(x_{2},t)-p_{vs}(x_{2},t)}{%
p_{i}(x_{2},t)}\right) $%
\end{tabular}%
$\ \medskip $

Le passage \`{a} la limite donne pour tout $x_{2}\in \lbrack 0;L]:%
\begin{tabular}{l}
$\underset{t\rightarrow +\infty }{\lim }\delta (x_{2},t):=\delta _{\infty
}(x_{2})$%
\end{tabular}%
\medskip $

o\`{u}%
\begin{equation}
\delta _{\infty }(x_{2})=\left\{
\begin{tabular}{l}
$1-\frac{\pi _{S}e^{\lambda _{i}(\frac{1}{\Theta _{S}}-\frac{1}{\Theta _{E}}%
)}}{p_{E}-2q\eta _{\Theta _{E}}x_{2}}$ si $x_{2}\leq \frac{p_{E}-p_{S}}{%
2q\eta _{\Theta _{E}}}$ \\
$1-\frac{\pi _{S}e^{\lambda _{i}(\frac{1}{\Theta _{S}}-\frac{1}{\Theta _{E}}%
)}}{p_{S}}$ sinon.%
\end{tabular}%
\right.  \label{cond}
\end{equation}

Remarque 4: \textit{Comme cette \'{e}tude s'interesse avant tout au
comportement asymptotique lorsque le temps }$t$ \textit{tend vers l'infini, }%
$t$\textit{\ est suppos\'{e} dans la suite assez grand pour pouvoir
approximer }$\delta (x_{2},t)$\textit{\ par }$\delta _{\infty }(x_{2}).$

\subsubsection{D\'{e}termination de la vitesse de recul de l'interface}

Pour un $z$ fix\'{e} dans $[0;1]$ et pour un instant $t$ fix\'{e}, soit $%
x_{2}=x_{2}(z,t)$ la position correspondante de l'interface entre le fluide r%
\'{e}siduel et le fluide condensable.

Suite \`{a} la condensation du fluide condensable pr\'{e}sent sur
l'horizontale de cote $z$ entre $x_{2}=0$ et $x_{2}=x_{2}(z,t)$ , cette
interface subit une recul\'{e}e dont la vitesse est $k_{0}(z,t)=%
\int_{0}^{x_{2}(z,t)}\delta (x_{2},t)dx_{2}$

\subsubsection{Comportement asymptotique de $x_{2}(t,z)$}

Comme dans le cas isotherme envisag\'{e} dans \cite{pous} et \cite{zel}, montrons que la
position de l'interface est born\'{e}e pour tout temps.

D'après \ref{vitesse} la vitesse est donnée par:
$$
{v}_{0}{(x}_{2},z,t)=w_{0}{(}z)-k_{0}(z,t),
$$
la position $x_{2}(z,t)$ \`{a} l'instant $t$ de l'interface entre les deux
fluides qui \'{e}tait en $M(x_{2}=0,z)$ \`{a} l'instant $t=0$ v\'{e}rifie
donc:\smallskip \medskip

\begin{tabular}{l}
$x_{2}(z,t)=w_{0}{(}z)t-\int_{0}^{t}k_{0}(z,\tau )d\tau =w_{0}{(}%
z)t-\int_{0}^{t}\int_{0}^{x_{2}(z,t)}\delta (x_{2},t)dx_{2}d\tau $%
\end{tabular}%
\medskip

Et puisque pour les grands temps nous avons:

$$
\delta (x_{2},t)\approx \delta _{\infty }(x_{2})\geq 1-\frac{\pi
_{S}e^{\lambda _{i}(\frac{1}{\Theta _{S}}-\frac{1}{\Theta _{E}})}}{p_{S}}>0
$$
la position $x_{2}(z,t)$ de l'interface est major\'{e}e par la solution
de l'\'{e}quation
$$
\xi _{2}=w_{0}{(}z)t-\left( 1-\frac{\pi _{S}e^{\lambda _{i}(\frac{1}{\Theta
_{S}}-\frac{1}{\Theta _{E}})}}{p_{S}}\right) \int_{0}^{t}\xi _{2}d\tau
$$

c'est-\`{a}-dire par:
$$
\frac{w_{0}{(}z)t}{1+\left( 1-\frac{\pi _{S}e^{\lambda _{i}(\frac{1}{\Theta
_{S}}-\frac{1}{\Theta _{E}})}}{p_{S}}\right) t}
$$

Et donc dans ce cas nous avons:
$$
\underset{t\rightarrow +\infty }{\lim }x_{2}(z,t)\leq \underset{%
t\rightarrow +\infty }{\lim }\frac{w_{0}{(}z)t}{1+\left( 1-\frac{\pi
_{S}e^{\lambda _{i}(\frac{1}{\Theta _{S}}-\frac{1}{\Theta _{E}})}}{p_{S}}%
\right) t}.
$$
soit finalement:
\begin{equation}
\underset{t\rightarrow +\infty }{\lim }x_{2}(z,t)\leq \frac{w_{0}{(}z)}{1-%
\frac{\pi _{S}e^{\lambda _{i}(\frac{1}{\Theta _{S}}-\frac{1}{\Theta _{E}})}}{%
p_{S}}}<+\infty.
\end{equation}
.

\subsubsection{\protect\bigskip Conclusion}

En se pla\c{c}ant dans des conditions  non
isothermes, avec un modèle très simple pour la température, le comportement global d\'{e}j\`{a} obtenu dans des conditions
isothermes dans  (cf. \cite{pous} et \cite{zel}) est maintenu:

S'il y a condensation du fluide entrant, la vitesse de l'interface entre le
fluide entrant et le fluide r\'{e}siduel d\'{e}croit progressivement vers $0$
et la position du front tend vers une position asymptotique finie quand le
temps $t$ \ tend vers l'infini.

\end{document}